\documentclass[conference]{IEEEtran}
\IEEEoverridecommandlockouts
\usepackage{amsmath,amsthm}
\usepackage{pgfplots}
\usepackage{amssymb}
\usepackage[fleqn,tbtags]{mathtools}
\usepackage{amsfonts}
\usepackage{graphicx}      % include this line if your document contains figures
\usepackage{xcolor}
\usepackage{hyperref}
\pgfplotsset{compat=1.16}
\def\rhopc{\rho_{\mathsf pc}}
\def\rhopci{\rho_{{\mathsf pc},i}}
\def\Upc{U_{\mathsf pc}}
\def\vpci{v_{{\mathsf pc},i}}
\def\Apci{A_{{\mathsf pc},i}}

\newtheorem{theorem}{Theorem}[section]

\newtheorem{remark}[theorem]{Remark}

\def\BibTeX{{\rm B\kern-.05em{\sc i\kern-.025em b}\kern-.08em
    T\kern-.1667em\lower.7ex\hbox{E}\kern-.125emX}}

\definecolor{matlabblue}{HTML}{0072BD}
\definecolor{matlaborange}{HTML}{D95319}
\definecolor{matlabyellow}{HTML}{EDB120}
\definecolor{matlabpurple}{HTML}{7E2F8E}
\definecolor{matlabgreen}{HTML}{77AC30}
\definecolor{matlablightblue}{HTML}{4DBEEE}
\definecolor{matlabred}{HTML}{A2142F}

\def\tt{^\intercal}

\hypersetup{
  colorlinks=true,
  urlcolor=matlabblue,% required dvipsnames optioon for xcolor
  citecolor=matlabgreen,% required dvipsnames optioon for xcolor
  linkcolor=matlabred% required dvipsnames optioon for xcolor
}

\begin{document}

%\title{Conference Paper Title*\\
%{\footnotesize \textsuperscript{*}Note: Sub-titles are not captured in Xplore and
%should not be used}
%\thanks{Identify applicable funding agency here. If none, delete this.}
%}

\title{
Low-order Linear Parameter Varying Approximations for Nonlinear
Controller Design for Flows
}

\author{\IEEEauthorblockN{1\textsuperscript{st} Amritam Das}
  \IEEEauthorblockA{\textit{Department of Electrical Engineering}\\
    \textit{Eindhoven University of Technology}\\
    P.O. Box 513, 5600 MB Eindhoven, The Netherlands \\
  \href{mailto:am.das@tue.nl}{am.das@tue.nl}
\href{https://orcid.org/0000-0001-8494-8509}{ORCID:0000-0001-8494-8509}}
\and
  \IEEEauthorblockN{2\textsuperscript{nd} Jan Heiland}
    \IEEEauthorblockA{%
      \textit{Max Planck Institute for Dynamics of Complex Technical Systems}\\
 39106 Magdeburg, Germany \\
  \href{mailto:heiland@mpi-magdeburg.mpg.de}{heiland@mpi-magdeburg.mpg.de}
\href{https://orcid.org/0000-0003-0228-8522}{ORCID:0000-0003-0228-8522}}
}

\maketitle

% \author{Amritam Das$^1$, Jan Heiland$^2$
% \thanks{$^1$Eindhoven University of Technology, Dept. of Electrical Engineering, Control Systems group, P.O. Box 513, 5600 MB Eindhoven, The Netherlands. E-mail: {\tt\small am.das@tue.nl}}
% \thanks{$^2$Max-Planck Institute for Dynamics of Complex Technical Systems, Sandtorstraße 1, 39106 Magdeburg, Germany. E-mail: {\tt\small heiland@mpi-magdeburg.mpg.de}}
% % \thanks{$^3$Vrije Universiteit Brussel (VUB), Dept. of Fundamental Electricity and Instrumentation, Pleinlaan 2, 1050 Brussels, Belgium}
% % \thanks{$^4$Eindhoven University of Technology, Dept. of Mechanical Engineering, Control Systems Technology group, PO Box 513, 5600 MB Eindhoven, The Netherlands}
% % \thanks{{\bf{Acknowledgements}:} This work has been carried out within the framework of the EUROfusion Consortium and has received funding from the Euratom research and training programme 2014-2018 under grant agreement No. 633053. The views and opinions expressed herein do not necessarily reflect those of the European Commission.}
% }

% \begin{document}

% \maketitle
% \thispagestyle{empty}
% \pagestyle{empty}

%%%%%%%%%%%%%%%%%%%%%%%%%%%%%%%%%%%%%%%%%%%%%%%%%%%%%%%%%%%%%%%%%%%%%%%%%%%%%%%%
\begin{abstract}
The control of nonlinear large-scale dynamical models such as the
incompressible Navier-Stokes equations is a challenging task.
The computational challenges in the controller design come from both the possibly large state space and the nonlinear dynamics. 
A general purpose approach certainly will resort to numerical linear algebra
techniques which can handle large system sizes or to model order reduction.
In this work we propose a two-folded model reduction approach tailored to
nonlinear controller design for incompressible Navier-Stokes equations and
similar PDE models that come with quadratic nonlinearities. Firstly, we
approximate the nonlinear model within in the class of LPV systems with a very
low dimension in the parametrization. Secondly, we reduce the system size to a
moderate number of states. This way, standard robust LPV theory for nonlinear
controller design becomes feasible. We illustrate the procedure and its
potentials by numerical simulations.
\end{abstract}


% For peer review papers, you can put extra information on the cover
% page as needed:
% \ifCLASSOPTIONpeerreview
% \begin{center} \bfseries EDICS Category: 3-BBND \end{center}
% \fi
%
% For peerreview papers, this IEEEtran command inserts a page break and
% creates the second title. It will be ignored for other modes.
\IEEEpeerreviewmaketitle

%%%%%%%%%%%%%%%%%%%%%%%%%%%%%%%%%%%%%%%%%%%%%%%%%%%%%%%%%%%%%%%%%%%%%%%%%%%%%%%%

\section{Introduction}

We consider general nonlinear, control-affine systems of type
\begin{equation}\label{eq:gen-cont-affine}
  \dot x(t) = f(x(t)) + Bu(t),
\end{equation}
where $x$ is the system's state with values in $\mathbb R^{n}$ with $n$ possibly large,
where $u$ is the input with values in $\mathbb R^{p}$, and where $f\colon
\mathbb R^{n}\to \mathbb R^{n}$. 

The computer-aided controller design for \eqref{eq:gen-cont-affine} with $n$ large is a
challenging problem with no generic approach being established yet.
Commonly used methods like \emph{backstepping} \cite{Kok92}, feedback
linearization \cite[Ch. 5.3]{Son98}, or sliding mode control
\cite{Dod15} require structural assumptions and, thus, may not be
accessible to a general computational framework. 
% A more ), or the repeated
% computation for subobtimal control laws like in model predictive control (MPC)
% \cite{GruP17}.
The both holistic and general approach via the \emph{HJB}
equations is only feasible for
very moderate system sizes or calls for model order reduction; see, e.g.,
\cite{BreKP19} for a relevant discussion and an application in fluid flow
control. 
As an alternative to reducing the systems size, one may consider approximations
to the solution of the \emph{HJB} of lower complexity. For that, e.g., truncated
polynomial expansions \cite{BreKP18} are considered or heuristic approximative solutions 
via the so called \emph{state-dependent Riccati equation}; see, e.g.,
\cite{BanLT07,HeiW23}.

In this work, we propose the embedding of \eqref{eq:gen-cont-affine} into the
class of \emph{linear parameter-varying} (LPV) systems of type
\begin{equation}\label{eq:gen-lpv-sys}
  \dot x(t) = A(\rho(t))\, x(t) + Bu(t),
\end{equation}
with $A(\rho(t))\in \mathbb R^{n\times n}$ for parameter values $\rho(t) \in
\mathbb R^{r}$ by implementing two layers of complexity reduction so that
established theory and algorithms (\cite{ApkGB95}) for robust LPV controller
design become available.

As laid out in \cite{HofW15}, the controller design techniques for {LPV} systems can
be classified into the categories \emph{polytopic}, \emph{{LFT}-based}, and
\emph{gridding}.

A general assumption of all these approaches is that the set of possible parameter
values is bounded. 
If, in addition, the parameter dependency of the coefficients is affine-linear,
then the theory and the related computations simplify significantly; see,
e.g., \cite{ApkGB95, Sch96}.

The polytopic approach is seen as the most developed approach with the notable
results from \cite{ApkGB95, BecP94} that provide algorithms and theory for a
scheduled $H_\infty$-robust controller and
that are the base of the \texttt{hinfgs} routine in the MATLAB \emph{Robust
Control Toolbox}. 

In a first approximation step, we seek for very low-dimensional encodings of the
states for replacing the nonlinear source terms by a low-dimensional linear
parameter varying (LPV) surrogate. 
In this step that is directed to adapt the nonlinearity to a format that is 
accessible to the LPV theory for computer-aided controller design, we tailor the
approximations for a best representation of the source terms at very
low-dimensions.

We will show that this approach can lead to LPV models that well approximate the
actual dynamics with as much as $r$ dimensions in the parametrization, with $r$
well less than $10$. 
Thus, application of standard linear matrix inequalities (LMI) approaches
already comes into reach. 
Still, this would require the solution of at least $r+1$ but more likely of
about $2^r$ coupled LMIs of the size of the original
system. 
Therefore, we propose a second layer of approximation that is tailored for the
low-order LPV representation to accurately follow the original dynamics
with a moderately sized state space.

% [TODO: Maybe some more words on the LPV controller design.]

This integrated approach to nonlinear controller design via tailored
approximations and established robust LPV theory has not been considered
so far. 

Direct relations to the vast research on LPV systems are given as follows.
Although the typically considered LPV systems are of moderate size (see \cite[Tab. V]{HofW15} that
classifies state space dimensions larger than $10$ as \emph{high dimensional}),
the unfavourable increase of complexity with the parameter dimension has
triggered various work on reducing the so-called \emph{scheduling dimension};
see, e.g., \cite{SirTW12}, \cite{KwiW08}, \cite{RizAV18}, or
\cite{KoeT20} that base on sparse optimization, \emph{principal component
analysis} (PCA), \emph{auto encoders}, or general \emph{deep neural networks}, respectively. 

The idea of using model order reduction in combination with LPV controller
design has been followed by \cite{HasW11} where a PDE system is reduced as a whole
before the treatment as an LPV system. Recently, we have delivered
a proof of concept (\cite{HeiW23}) that LPV approximations of nonlinear systems
can with low parameter dimensions work well for nonlinear controller design.

\section{POD for low-dimensional LPV approximations}

% For this first prove of concept we will consider the established model reduction
% technique of \emph{Proper Orthogonal Decomposition} (POD). 
Under mild conditions (see, e.g., \cite{BenH18}), the system
\eqref{eq:gen-cont-affine} can be
  brought into the so-called \emph{state-dependent coefficient} (SDC) form
\begin{equation}\label{eq:SDC-quad-system}
  \dot x(t) = A(x(t))\, x(t) + Bu(t),
\end{equation}
with $A(x(t)) \in \mathbb R^{n\times n}$.
We will assume that, in particular, the map $x \to A(x) \in \mathbb R^{n\times
n}$ is (affine)
linear, i.e. $A(x)$ can be realized as $A(x) = A_0 + L(x)$ with $L$ linear.

\begin{remark}
  For some systems like
  the spatially-discretized Navier-Stokes equations the SDC representation
  with affine dependencies is naturally induced by its structure; see
  \cite{HeiBB22}.
  If the map $x\to A(x)$ is not linear, one can seek to find an approximation
  that is linear; see, e.g., \cite{KoeT20}.
\end{remark}

We note that \eqref{eq:SDC-quad-system} is an LPV representation
\eqref{eq:gen-lpv-sys} of \eqref{eq:gen-cont-affine} with the trivial
parametrization $\rho(t)=x(t)$ and with, in particular, a large parameter
dimension namely $r=n$.

Next we illustrate how a general model order reduction scheme, can provide LPV
approximations with possibly low-dimensional, e.g., $r\ll n$, parameter domains.

Let $V_r \in \mathbb R^{n\times r}$ be a POD basis that encodes and decodes the
state $x(t)$ as $\rho(x(t)) = V_r^Tx(t)$ and 
\begin{equation*}
  x(t) \approx \tilde x(t) = V_r\rho(x(t)) = \sum_{i=1}^r v_i \rho_i(x(t)),
\end{equation*}
where $v_i\in \mathbb R^{n}$ is the $i$-th POD mode.

With this encoding and decoding, the nonlinear term in system
\eqref{eq:SDC-quad-system} can be approximated as
\begin{equation*}
  \begin{split}
  A(x(t))\,x(t) \approx A(\tilde x(t))\,x(t) = \\ 
  A(\sum_{i=1}^r v_i \rho_i(t))\,x(t)
  = [A_0 + \sum_{i=1}^r \rho_i(t) L(v_i)]\, x(t) 
  \end{split}
\end{equation*}
which defines an affine-linear low-dimensional LPV approximation
\begin{equation}\label{eq:ldlpv-system}
  \tilde x(t) = [A_0 + \sum_{i=1}^r \rho_i(\tilde x(t)) L(v_i)]\, \tilde x(t) +
  Bu(t)
\end{equation}
to \eqref{eq:SDC-quad-system} and \eqref{eq:gen-cont-affine}.
We set
\begin{equation}
  A_i := L(v_i), \quad i = 1, \dotsc, r.
\end{equation}

In the second step, we now project the LPV approximation to reduced order
coordinates.
For that let $V_k$ with $k\geq r$ be the POD basis that includes $V_r$ and $\bar
\rho(t) = V_k \tilde x(t)$. 
Then an approximation to \eqref{eq:ldlpv-system} with state dimension $k$ (as
opposed to $n$) reads
% , additionally, the system itself is projected to the system of $k\geq r$ POD
% coordinates $\bar \rho = V_k x$, an LPV-POD approximation is given as
\begin{equation}\label{eq:pod-ldlpv-sys}
  \dot {\bar \rho}(t) = [\bar A_0 + \sum_{i=1}^r \bar \rho_i(t) \bar A_i]\, \bar
  \rho(t) + \bar Bu(t).
\end{equation}\\\
with $\bar A_i := V_k^TA_iV_k$, for $i=0,1,\dots,r$ and $\bar B:= V_k^TB$.

Note that $r$ can be chosen independently of $k$. Thus, a very low-dimensional approximation
can be achieved independently of a possibly larger state space that can be
tailored for the best compromise in terms of size and accuracy.

\section{Controller Design for LPV Systems}

The polytopic approach is seen as the most developed approach with the notable
results from \cite{ApkGB95, BecP94} that provide algorithms and theory for a
scheduled robust controller and
that are the base of the \texttt{hinfgs} routine in the MATLAB \emph{Robust
Control Toolbox}. 

\section{Implementation Issues}

For the computation of the snapshots that are used to extract the POD basis, a
forward simulation for an example input is performed.

\subsection{Polytope or Bounding Box}
The computed snapshots are then projected to the $\rho$-coordinates in order to
estimate the polytope $W\subset \mathbb R^{r}$ that contains $\rho(t)$. 

The vertices $w_i$ of $W$ are then used to define the controller.

Intuitively, the performance of the controller will be better if the volume
$|W|$ of $W$ is smaller and if the vortices $w_i$ are closer to extremal values
of $\rho(t)$. 
In this respect, the optimal choice of $W$ would be the convex hull of the
snapshots of $\rho$.

On the other hand, the convex hull is likely to have a large number of vortices,
which makes the application of LPV controller design approaches costly as they,
e.g., require the solution of $n$ coupled LMIs where $n$ is the number of
vertices of $W$. 

Typically, $W$ is chosen as a bounding box with, accordingly, $2^r$ vortices. In
order to reduce the volume, the box can be rotated and expressed in coordinates
obtained of the principal components; see, e.g., \cite{KwiW08}. 
In our affine-linear case \eqref{eq:ldlpv-system}
with
\begin{equation}
  \rho (t) = \Upc\Upc^T \rho(t)=: \Upc \rhopc(t) 
\end{equation}
and, accordingly, 
\begin{equation}
  \tilde x (t) = V_r \rho(t)= V_r \Upc \rhopc(t) =: \sum_{i=1}^r \vpci
  \rhopci(t)
\end{equation}
where $\Upc$ is the \emph{principal components} coordinate transformation,
this reparametrization reads
\begin{equation}
 \sum_{i=1}^r \rho_i(t) A_i = L(\tilde x(t)) = \sum_{i=1}^r \rhopci(t) \Apci ,
\end{equation}
where $\Apci:=L(\vpci)$ and $\vpci$ denoting the $i$-th column of $V_r\Upc$, for
$i=1,\dotsc,r$. For a direct retransformation of the system, we can also resort to the relation
\begin{equation}
  \Apci = \sum_{j=1}^r U_{ji}A_j,
\end{equation}
where $U_{ji}$ is the $j$-th row entry of the $i$-th column of $\Upc$.

Adding on these standard ways on defining the parameter domains, in our examples we use optimization to find a suitable polytope that gives a
good compromise of volume and number of vortices and that underbids the bounding
box in both variables. For this, the following optimization setup was defined
and solved using the built-in methods of computing convex hulls and genetic
optimization in \emph{SciPy}; see, e.g., \cite{scipy2020}.

\begin{itemize}
  \item[(0.)] Let $V\in \mathbb R^r$ be (the set of vertices of) the convex hull
    of given measurements of $\rho(x(t_j))$, for
    $j=1,\dots, N$ and a given state trajectory $x$. 
  \item[(i.)] In the $i$-th iteration, extend $V$ by $n_k$ vertices $v^{(i)}_k\in \mathbb R^{r}$ and compute
    the convex hull $V^{(i)}$ of $V \cup \{v^{(i)}_1, \dotsc, v^{(i)}_{n_k}\}$
  \item[(ii.)] Update $v^{(i)}_k$ to minimize both the volume and
    the number of vertices of $V^{(i)}$.
\end{itemize}

\begin{figure}
  \includegraphics[width=\linewidth]{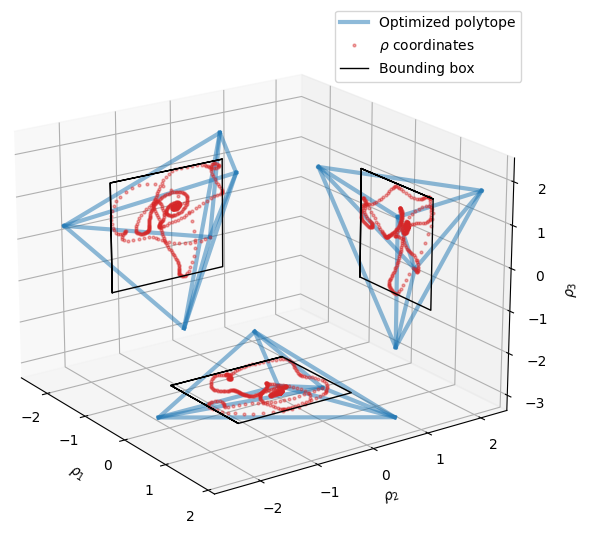}
  \caption{An illustration of the first three components of $\rho(t)$, of the bounding box (with 8 vortices), and an optimized
polytope of less volume and with only 5 vortices). Each data object is in
three-dimensional space but projected along the coordinate axes to planes in the
two-dimensional space.}
  \label{fig:3d-illustration-bb-polytope}
\end{figure}

We report on efficiency (and feasibility) of the computation of the LPV
controller and on its performance for the three approaches of considering
\begin{itemize}
  \item the bounding box in the original $\rho(t)$ coordinates,
  \item the bounding box in the PCA coordinates of  $\rho(t)$, or
  \item an optimized polytope of fewer vertices 
\end{itemize}
as the polytope for the parameter variation; see Figure
\ref{fig:3d-illustration-bb-polytope} for an illustration of the bounding box
and an optimized polytope.

\section{Numerical Example}

\begin{figure}
  \includegraphics[width=\linewidth]{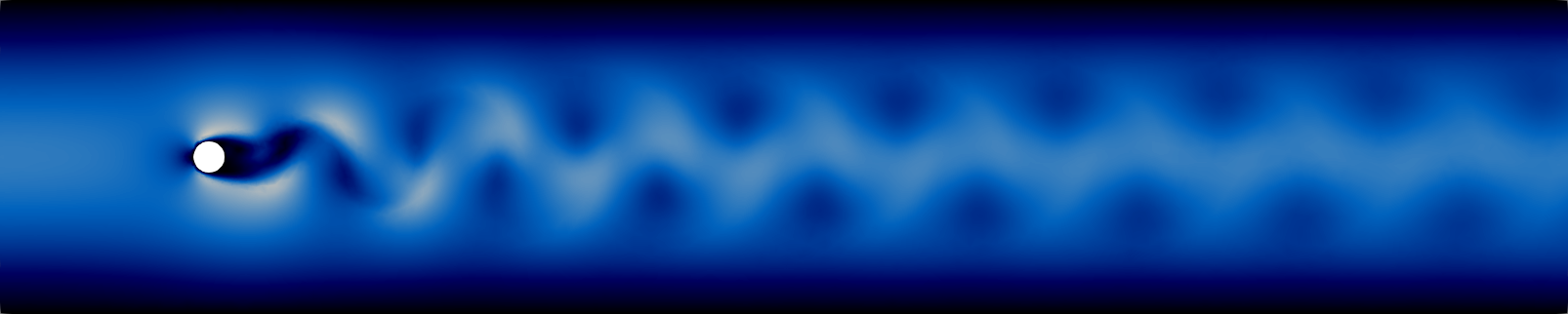}
  \caption{Illustration of the computational domain and the developed state of
  the uncontrolled flow field.}
  \label{fig:cw-domain-flow}
\end{figure}

We consider the two-dimensional cylinder wake with control at moderate Reynolds
numbers; see Figure \ref{fig:cw-domain-flow}. 
The controls are designed as two outlets at the cylinder periphery of size
$\pi/6$ located at $\pm \pi/3$ through which
fluid can be injected or sucked away. 
Mathematically, the control is modelled by parabola-shaped spatial shaped
function that is scaled by the scalar control value.
For inclusion in the FEM scheme, these Dirichlet conditions are relaxed towards
Robin-type boundary conditions with a relaxation parameter $\epsilon=10^{-5}$; see
\cite{BehBH17} for implementation details.

As the output $y$, we consider averaged velocities in 3 square domains of
observations of size $D^2$ located at a distance of $H$ behind the cylinder
symmetrically with respect to the channel middle. Here $D$ denotes the diameter
of the cylinder and $H$ the height of the channel. With the two components of
the velocity, overall, an output $y(t)\in \mathbb R^{6}$ is obtained with the first
3 values corresponding to the stream wise components and the final 3 values to
the lateral components of the velocities.

The corresponding PDE model is spatially discretized by \emph{Taylor-Hood}
quadratic-linear mixed finite elements on a nonuniform grid. 
From the FEM model of about 50\,000 degrees of freedom, the low-dimensional LPV
approximation is obtained through the following algorithm

\begin{enumerate}
  \item Starting from the steady-state solution, the FEM model is integrated in
    time from $t=0$ to $t=5$ with a test input applied to trigger the
    instabilities.
  \item From the FEM solution, 417 equispaced snapshots are collected to define
    the POD basis.
  \item With the POD basis, the reduced order LPV approximation as in
    \eqref{eq:pod-ldlpv-sys} is computed.
\end{enumerate}

\begin{figure}
  \includegraphics[width=\linewidth]{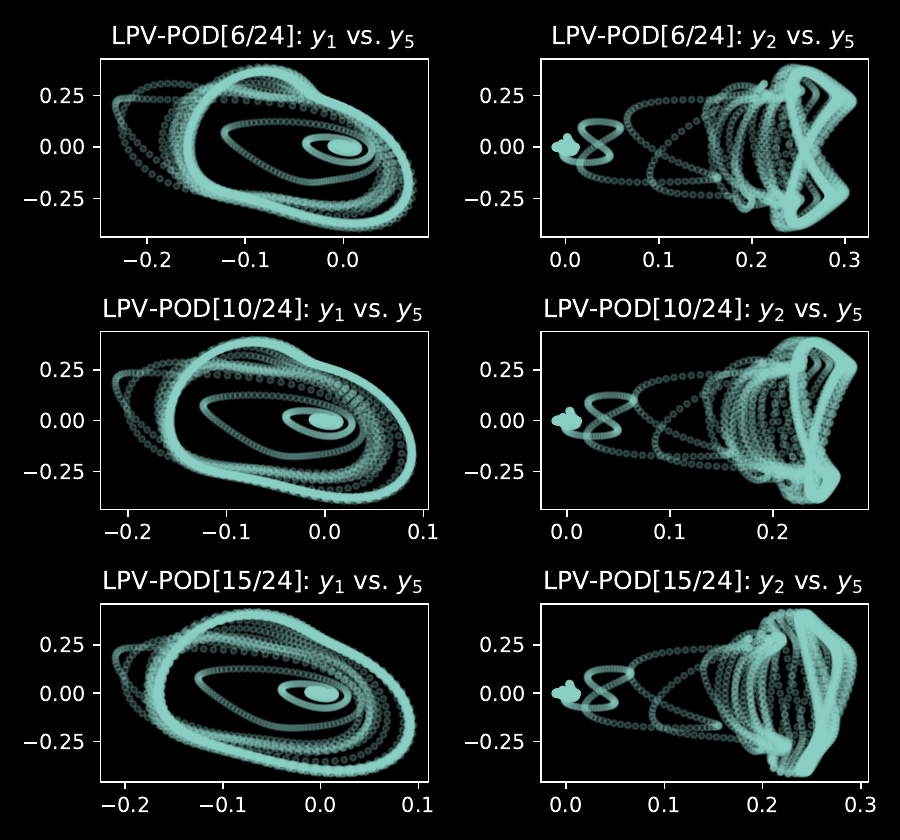}
  \caption{Phase portraits $y_2$ vs. $y_5$ and $y_1$ vs. $y_5$ of the output $y$ for POD
  dimension $k=24$ and parameter dimension $r=6, 10, 15$ on the time
interval $(0,12)$.}
  \label{fig:pps-k24}
\end{figure}

\begin{figure}
  \includegraphics[width=\linewidth]{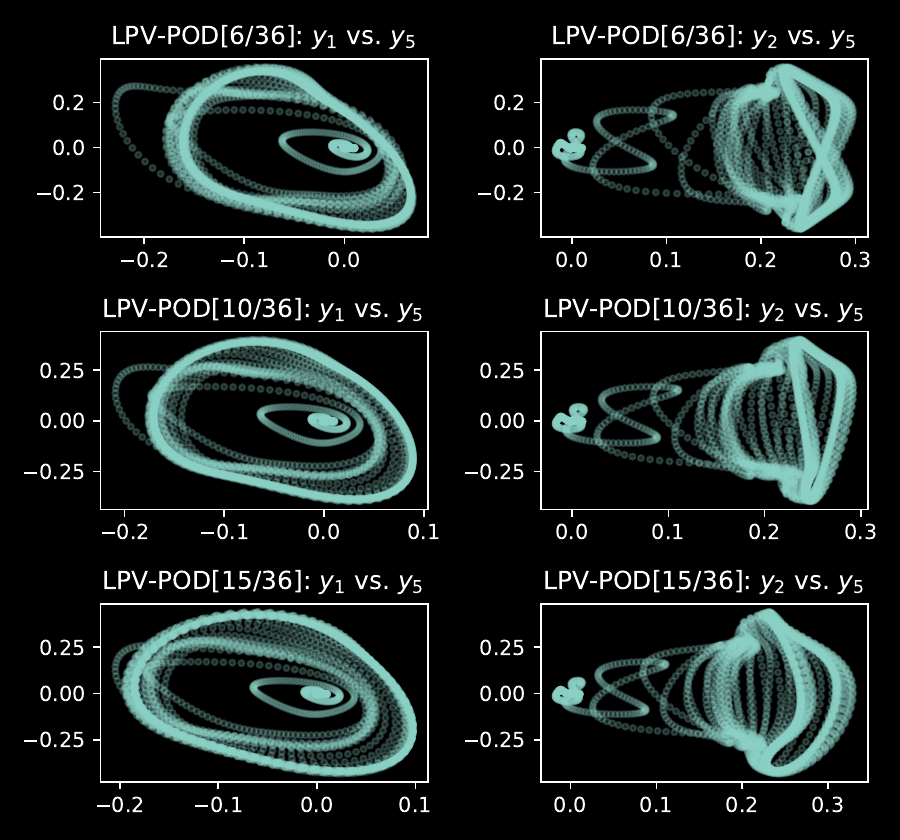}
  \caption{Phase portraits $y_2$ vs. $y_5$ and $y_1$ vs. $y_5$
  of the output $y$ for POD dimension $k=36$ and parameter dimension $r=6, 10,
15$ on the time interval $(0,12)$.}
  \label{fig:pps-k36}
\end{figure}

\begin{figure}
\includegraphics[width=\linewidth]{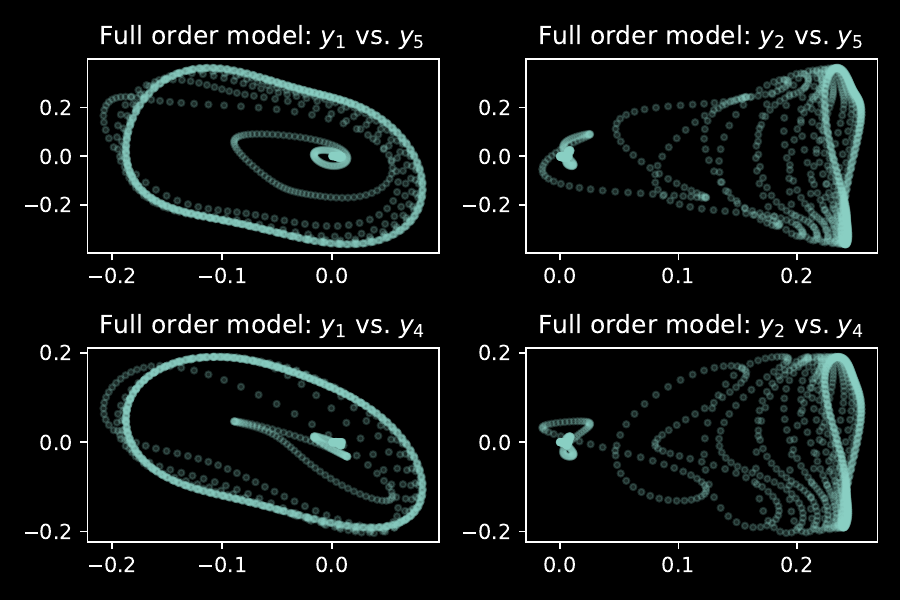}
  \caption{Phase portrait of selected components of the output $y$ for the full
  order FEM model on the time interval $(0, 12)$. Because of the symmetry of the
  overall setup, the phase portraits of the remaining
components bear no additional information and are not shown here.}
  \label{fig:pps-fom}
\end{figure}

In view of determining $k$ and $r$, i.e., the size of the POD reduced model and
the size of the parameter in the LPV approximation, we note that the system is
chaotic which makes a quantitative decision delicate (as small perturbations
have arbitrarily large effects). 
That's why we took the qualitative view of examining the resulting limit cycles
of selected components in the phase portraits of $y_2$ vs. $y_5$ and $y_1$ vs.
$y_5$, see Figures \ref{fig:pps-k24}--\ref{fig:pps-k36} for results of different
choices of $k$ and $r$ and Figure \ref{fig:pps-fom} for the reference. 
These phase portraits show the data points, say $(y_2(t_i), y_5(t_i))$, for
$500$ time equidistant instances $t_i$ from the output $y$ of the corresponding
simulations with zero inputs on the time interval $(0, 12)$.

For the following numerical studies we chose the setup $k=36$ and $r=6$ that,
judging from Figure \ref{fig:pps-k36} (first line) in comparison to Figure
\ref{fig:pps-fom} (first line), seems to well cover at least the range of
values in the phase portraits for the smaller values of $r$.

% we investigated the approximation quality of
% the POD reduced order model for several choices of $k$ on the training interval
% $(0,1)$ and, subsequently, identified $r$ through examining the residual
% \begin{equation}
%   \int_0^5  \bigl \|\sum_{i=1}^r\rho_i(t) A_i \rho(t) - L(\rho(t))\rho(t)\bigr \| dt
% \end{equation}
% with $\rho$ as the solution of the POD reduced order model on $(0, 5)$.

\subsection{Computing the LPV Controller}

For the computation of the LPV controller we used the \emph{Matlab Robust
Control Toolbox}\cite{MatRC} in the release 2022b with built-in function \texttt{hinfgs} that computes an LPV controller with a
guaranteed quadratic robustness performance $\gamma_*$.
The computational costs of the underlying optimization with coupled LMI
constraints are significant and make larger values of $k$ and $r$ quickly
infeasible for numerical studies.

As for the different approaches of defining the enclosing polytope
$W\supseteq\{\rho(t)\colon t>0\}$ we made the following observations.

\begin{figure}
  \includegraphics[width=\linewidth]{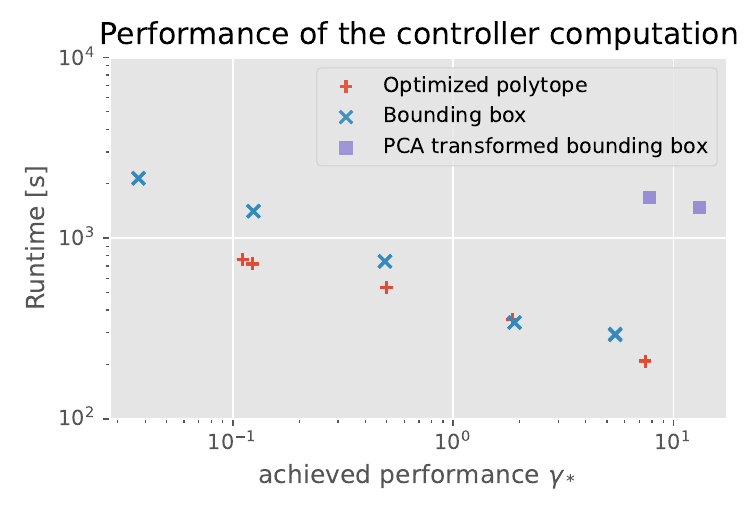}
  \caption{Achieved $\gamma_*$ versus CPU time (values of three consecutive simulations)
    for the bounding box, the PCA transformed coordinates, and the optimized
    vertices setup. The left most values correspond to the best achieved
  $\gamma_*$ for which the iteration stagnates.}
  \label{fig:cpuvsgammas}
\end{figure}

\begin{itemize}
  \item The bounding box in the original coordinates with $2^6=64$ achieves the
    best robustness performances $\gamma_*$ though with rather high
    computational costs.
  \item The bounding box in the PCA coordinates comes with the same number of
    vertices to consider.
    % and, thus, with the same computational costs per iteration. 
    However, the observed convergence in $\gamma_*$ was tediously slow which led
    to infeasible numerical costs for lower target values of $\gamma_*$.
  \item Using an optimized polytope of $20$ vertices, each iteration in the
    \texttt{hinfgs} computation was sped up by a factor of $3$ which well
    compensated for an overall slower convergence. 
    However the slower convergence even led to stagnation so that the best
    achievable values of $\gamma_*$ were larger than that of the bounding box
    approach.
\end{itemize}
The results on the performance of \texttt{hinfgs} for the different setups and
for different levels of $\gamma_*$ are displayed in the chart of Fig.
\ref{fig:cpuvsgammas}.
The conducted experiments suggest that a compact representation of $W$ (in the
sense that its vertices are evenly distributed in space or that the ratio of
surface over volume of $W$ is rather small) is beneficial for
convergence in $\gamma_*$. 
Apparently, the vertices of the optimized polytope are further apart (see Fig. 
\ref{fig:3d-illustration-bb-polytope}) which may explain the slower convergence
and earlier stagnation. 
An explanation for the minor performance of the PCA coordinate transformation
may call on the interpretation of the PCA concentrating the variance of the data
in the leading principal components. This may result in vertices that widely
distributed in one dimension and closely located in another. Nonetheless, the
strong (in this case negative) effect supports the
idea that a transformation of the $\rho$ coordinates has the potential of improving the
performance of \texttt{hinfgs}.

\subsection{Nonlinear Closed Loop Simulations}

Since the robust control toolbox in \emph{Matlab} does not provide the
functionality to evaluate the LPV controller within a general polytope%
\footnote{Basically, \emph{Matlab} has no built-in function to compute
\emph{barycentric coordinates} in a convex polytope in dimensions higher than
$r=3$.}%
, we considered closed loop simulations with the controller obtained through the
bounding box approach. 
Also the built-in simulation routines only support predefined parameter
trajectories so that the closed-loop system was set up manually and simulated
with the time integrator \texttt{ode15s}.

\begin{figure}
  \includegraphics[width=\linewidth]{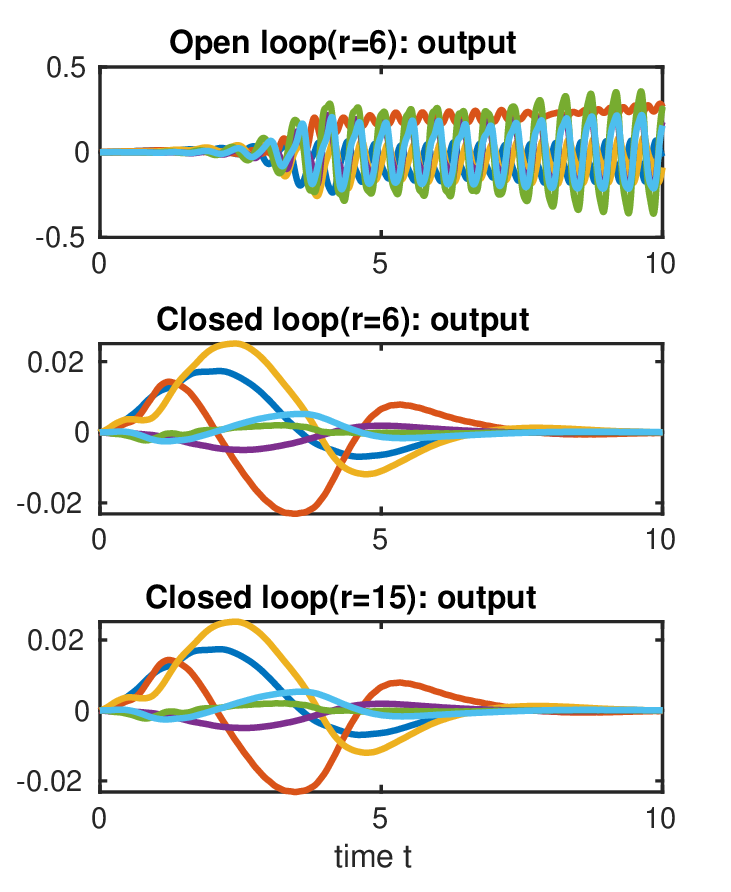}
  \caption{The open loop output of the model for $r=6$ and the closed loop for
  LPV controller designed for the $r=6$ model applied to this model and used in
the larger ($r=15$) model.}
\label{fig:olcl-outputs}
\end{figure}

As the result, this nonlinear controller did well stabilize the nonlinear system (with $r=6$
that it was built upon) as illustrated
in Fig. \ref{fig:olcl-outputs}. Also, this controller proved a certain
robustness by performing similarly well for the model with parameter dimension $r=15$ which is
similar but also shows different dynamical patterns; cp. Fig.
\ref{fig:pps-k36}(first row vs. third row).

\subsection{Remarks on Nonlinear Systems as LPV Systems}

Although the embedding of a nonlinear system via $A(\rho(t)) = A(\rho(x(t))$
into the class of LPV systems is readily covered by the available LPV theory, it
comes with practical consequences in the controller definition. 
Firstly, as the trajectory of $\rho$ is not preset but defined through the state
trajectory $x$, the containing polytope (or bounding box) has to be estimated
from state estimations or approximations. 

Secondly, the effect of the feedback is ambivalent. While a functioning
controller will stabilize the system around the working point and, thus, prevent
$x$ and $\rho(x)$ from attaining extreme values, short-term perturbations may
lead to overshoots which can drive the system out of the estimated region.

While, practically, such an overshoot can often be compensated, the computed
controller will fail immediately as the parameter update for the controller is
no more well-defined. In our experiments we observed these critical overshoots
when considering nonzero initial values or discontinuous disturbance signals. 
Therefore, in the presented simulations, we started from the zero initial conditions and applied a
    disturbance that smoothly faded out after $t=2$.

\subsection{Code Availability}
The LPV system data (ready for import in Matlab) and the scripts that were
used to obtain the presented numerical results are available for immediate
reproduction from 
\href{https://dx.doi.org/10.5281/zenodo.10073483}{doi:10.5281/zenodo.10073483}
under a CC-BY license.

\section{Conclusions and Outlook}

In this work, we have employed a two-level model order reduction approach so
that, eventually, established LPV theory and algorithms become available for
nonlinear controller design for general nonlinear systems like the
incompressible Navier-Stokes equations.

In a numerical experiment, we illustrated the potential and feasibility of the
combined approach and identified pitfalls of the approach as well as limits in
the existing functionality of standard control systems software. 

Notably, although the theory applies and although favourable properties of
tailored polytopes have been illustrated, routines for LPV
controller synthesis basically only allow for bounding boxes as enclosing
polytopes. Thus, a future work will concern interpolation of controllers in
polytopes using, e.g., the formulas provided in \cite{Wac11}.

Another immediate future development could concern the solution of large-scale
linear matrix inequalities in the context of LPV controller design.

\bibliographystyle{IEEEtran}
% \bibliography{bibfile}
\bibliography{ecc-nse-lpv}
\end{document}